\documentclass[11pt,amssymb]{article}

\usepackage{epsf}

\usepackage{amssymb,amsthm,verbatim}

  \font\tenmsb=msbm10

\font\sevenmsb=msbm7 \font\fivemsb=msbm5

\newfam\msbfam

\textfont\msbfam=\tenmsb

\scriptfont\msbfam=\sevenmsb

\scriptscriptfont\msbfam=\fivemsb

\def\Bbb#1{\fam\msbfam\relax#1}

\newtheorem{thm}{Theorem}[section]

\newcommand{\ep}{\epsilon}
\newtheorem{prop}[thm]{Proposition}

\newtheorem{cor}[thm]{Corollary}

\newtheorem{lem}[thm]{Lemma}

\newtheorem{conj}[thm]{Conjecture}

\newtheorem{exa}[thm]{Example}

\newtheorem{defn}[thm]{Definition}

\newtheorem{clm}[thm]{Claim}

\newtheorem{rem}[thm]{Remark}

\newtheorem{obs}[thm]{Observation}

\newtheorem{qs}[thm]{Quesion}

\newcommand{\ben}{\begin{enumerate}}

\newcommand{\een}{\end{enumerate}}

\newcommand{\brem}{\begin{rem}}

\newcommand{\erem}{\end{rem}}

\newcommand{\blem}{\begin{lem}}

\newcommand{\elem}{\end{lem}}

\newcommand{\bcl}{\begin{clm}}

\newcommand{\ecl}{\end{clm}}

\newcommand{\bthm}{\begin{thm}}

\newcommand{\ethm}{\end{thm}}

\newcommand{\bq}{\begin{qs}}
\newcommand{\eq}{\end{qs}}

\newcommand{\bpr}{\begin{prop}}

\newcommand{\epr}{\end{prop}}

\newcommand{\bco}{\begin{cor}}

\newcommand{\eco}{\end{cor}}

\newcommand{\bcon}{\begin{conj}}

\newcommand{\econ}{\end{conj}}

\newcommand{\bde}{\begin{defn}}

\newcommand{\ede}{\end{defn}}

\newcommand{\bex}{\begin{exa}}

\newcommand{\eexa}{\end{exa}}

\newcommand{\bexe}{\begin{exe}}

\newcommand{\eexe}{\end{exe}}

\newcommand{\bobs}{\begin{obs}}

\newcommand{\eobs}{\end{obs}}

\begin{document}

\title{Uniform Kazhdan Constant for some families of linear groups}  
\author{Uzy Hadad}
\date{\today}                   
\maketitle                      


\begin{abstract} Let $R$ be a ring generated by $l$ elements with stable range $r$. Assume that the group $EL_d(R)$
has Kazhdan constant $\ep_0>0$ for some $d \geq r+1$. We prove that
there exist $\ep(\ep_0,l) >0$ and $k \in \mathbb{N}$, s.t. for every
$n \geq d$, $EL_n(R)$ has a generating set of order $k$ and a
Kazhdan constant larger than $\epsilon$. As a consequence, we obtain
for $SL_n(\mathbb{Z})$ where $n \geq 3$, a Kazhdan constant which is
independent of $n$ w.r.t generating set of a fixed size.
\end{abstract}

\section{Introduction}

Let $k\in N$ and $ 0 < \ep \in \mathbb{R}$. A group $\Gamma$ is
said to have Kazhdan property $(T)$ with Kazhdan constant
$(k,\ep)$, if $\Gamma$ has a set of generators $S$, with $|S| \leq
k$ satisfying:

If $(\rho, \mathcal{H})$ is a unitary representation of $\Gamma$,
$\rho:\Gamma\rightarrow \mathcal{U(H)}$, with unit vector $v\in
\mathcal{H}$, such that for all $s \in S$, $\|\rho(s)v-v\| \leq
\ep$, then $\mathcal{H}$ contains a non-zero $\Gamma$-invariant
vector.

In response to a question raised by Serre, Shalom \cite{Sh1} and
Kassabov \cite{Kas2} showed that if one takes the elementary
matrices as the set of generators, then one gets that
$(2(n^2-n),\ep_n)$ with $(42 \sqrt{n}+860)^{-1} \leq \ep_n <
2n^{-\frac{1}{2}}$, is a Kazhdan constant for $SL_n(\mathbb{Z})$.
Here we show that by taking a different set of generators, $k$ and
$\ep$ can be made independent of $n$:

\bthm \label {main_z} There exist $k \in N$ and $0<\epsilon \in
\mathbb{R}$ s.t. for every $n\geq 3$, $SL_n(\mathbb{Z})$ has
Kazhdan constant $(k,\ep)$.
 \ethm

Theorem \ref{main_z} is deduced from a much more general result:

\bthm  \label{main_them0} Let $R$ be an associative ring generated
by $a_1,...,a_l$ with stable range $r$, and assume that for some
$d \geq r+1$, the group $EL_d(R)$ has Kazhdan constant
$(k_0,\ep_0)$. Then there exist $\ep=\ep(\ep_0,l)> 0$ and
$k=k(k_0,l) \in \mathbb{N}$, s.t. for every $n \geq d$, $EL_n(R)$
has Kazhdan constant $(k,\ep)$.

\ethm

Although we require $ d \geq r+1$, the result is true also if
$d<r+1$ provided that $EL_m(R)$ is a bounded product of
$EL_{m-1}(R)$ for every $r+1 > m
> d$.

For the definition of stable range and $EL_d(R)$ see Section
\ref{back}. As $EL_3(\mathbb{Z})=SL_3(\mathbb{Z})$ has property
$(T)$ and the stable range of $\mathbb{Z}$ is $2$, Theorem
\ref{main_z} follows from Theorem \ref{main_them0}.

Shalom \cite{Sh3,Sh4} proved that for $n > l+2$,
$SL_n(\mathbb{Z}[x_1,...,x_l])$ has Kazhdan property $T$ (see
Theorem \ref{Shalom} in this paper). Therefore we get:

\bco \label {Z_x_them} There exist $k=k(l)$ and $\ep=\ep(l)$ such
that for every $ n > l+2$, $SL_n(\mathbb{Z}[x_1,...,x_l])$ has
Kazhdan constant $(k,\ep)$. \eco

Since the stable range of any ring of integers $\mathcal{O}$ in a
global field is $2$, and it is known that $SL_3(\mathcal{O})$ has
property $T$ (see \cite{Sh1}), a similar result holds:

\bco \label {main_them1} For every ring of integers $\mathcal{O}$
in  a global field $K$, there exist $k=k(\mathcal{O})$ and
$\ep=\ep(\mathcal{O})$ such that for every $ n \geq 3$,
$SL_n(\mathcal{O})$ has Kazhdan constant $(k,\ep)$.  \eco

\brem If $\mathcal{O}$ is generated by $l$ elements, then
$SL_n(\mathcal{O})$ is a quotient of
$SL_n(\mathbb{Z}[x_1,...,x_l])$. Therefore for $n > l+2$,
$k(\mathcal{O})$ and $\ep(\mathcal{O})$ depend only on the number
$g(\mathcal{O})$ of generators of $\mathcal{O}$ as a ring. Also it
is known that $g(\mathcal{O})$ depends only on the discriminant of
$\mathcal{O}$ (see \cite{Pl}).\erem

The idea of this work is inspired by the work of Shalom in
\cite{Sh1,Sh2,Sh3} who relates property $(T)$ to bounded generation
and stable range, and also from the work of Kassabov in \cite{Kas1}
who extended Shalom's results and proved that there exist $k \in
\mathbb{N}$ and $ \ep>0$ s.t. for any finite commutative ring $R$,
and any $n \geq 3$, the group $EL_3(M_n(R))$ has Kazhdan constant
$(k,\ep)$ independent of $n$. Their proof is based on the fact that
the group in question is boundedly generated by elementary matrices (see \cite{CK} for the bounded elementary generation property of the group $SL_n(\mathcal{O})$ for $n \geq 3$).
In our case we show that the groups in question are boundedly
generated by elementary matrices and a common group with Kazhdan
property $(T)$.

\section{Background} \label{back}

\subsection{Property (T)}

 Property (T) was introduced by Kazhdan in \cite{Kaz}.
Since then, it has found numerous applications in various areas of
mathematics. Among them, for example, are constructions of
expander graphs \cite{Mar} , the product replacement algorithm \cite{LP} and bounds on mixing time of random walks on groups (see [Ch], [L1] for references and details).

 \bde Let
$\Gamma$ be a discrete group, $S \subset \Gamma $ a subset, $\ep
>0$, and $(\rho,\mathcal{H})$ be a unitary representation of the
group $\Gamma$. A vector $0\neq v \in \mathcal{H}$ is called
$(S,\epsilon)$-invariant, if $\|\rho(g)v-v\|\leq \epsilon \|v\|$
$\forall g\in S$. A discrete group $\Gamma$ is said to have Kazhdan
property $(T)$, if there exist a finite set $S \subset \Gamma$ and $
\ep >0$, such that every unitary representation with
$(S,\epsilon)$-invariant vector, contains a non-zero
$\Gamma$-invariant vector. In that case $(|S|,\ep)$ is called a
Kazhdan constant for $\Gamma$. We also sometimes say that $\Gamma$
has Kazhdan constant $\ep$ w.r.t the generating set $S$.\ede

We will need the following known Lemma (see Lemma \textbf{2.2} in
\cite{Sh3} and also Lemma \textbf{1} in \cite{KLN}).

 \blem \label{BP_gp}
If $\Gamma$ is a product of subgroups $H_1,H_2,...,H_k$, i.e.
$\Gamma=H_1 \cdot H_2 \cdot \cdot \cdot H_k$, where each $H_i$ has
Kazhdan constant $\ep_0$ w.r.t the generating set $S_i$ then
$\Gamma$ has a Kazhdan constant $\ep=\ep(\ep_0,k)$ w.r.t their union
$S=\bigcup S_i$. \elem

For more information and introduction for the subject we refer the
reader to ~\cite{L1}.

\subsection{$EL_n(R)$ and Property $(T)$} \label{EL_T}

\textbf{Notations:}
 Let $R$ be an associative ring with unit which is
generated by the elements $a_1,...,a_l$. Let $r \in R$ and let $i,j
\in \mathbb{N}$ s.t. $1 \leq i \neq j \leq n$. Denote by $e_{ij}(r)$
the $n \times n$ matrix with $1$ along the diagonal, $r$ in the
$(i,j)$ position, and zero elsewhere. Note that $e_{ij}(-r)$ is the
inverse of $e_{ij}(r)$ so that $e_{ij}(r) \in GL_n(R)$. These are
the \emph{elementary matrices}. The subgroup of $GL_n(R)$ which they
generate is the elementary group $EL_d(R)$.

The group $EL_d(R)$, provided  $d \geq 3$, is generated by the set
$S_d(R)$, where $S_d(R)$ contains the set of $2(d^2-d)$ elementary
matrices with $\pm 1$ off the diagonal and the set of $4l(d-1)$
elementary matrices $e_{ij}(\pm a_m)$ with $|i-j|=1$ and $1\leq
m\leq l$.

\bde The group $\Gamma=EL_d(R)$ is said to have the bounded
elementary generation property if there is a number $N=BE_d(R)$
such that every element of $\Gamma$ can be written as a product of
at most $N$ elementary matrices. \ede

 In [Kas1, Theorem \textbf{5}] Kassabov generalized a result of Y. Shalom \cite{Sh1} and proved the following:

\bthm \label {kassabov} Suppose that $d \geq 3$ and $R$ is a
finitely generated associative ring such that $EL_d(R)$ has the
bounded elementary generation property. Then $EL_d(R)$ has
property $T$ with an explicit lower bound for the Kazhdan constant
of $EL_d(R)$ with respect to the generating set $S_d(R)$.

\ethm

Moreover, Kassabov, in his proof of the above theorem, proved the
following lemma (Lemma 1.1 in \cite{Kas1} see also Corollary 3.5 in \cite{Sh1}) which plays a crucial
roll in this paper.

\blem \label{Kas_el} There exists a constant
$M(l)<3\sqrt{2}(\sqrt{l}+3)$ such that every f.g associative ring
$R$, which is generated by $l$ elements and for every $k \geq 3$,
the group $EL_k(R)$ satisfies the following property: Let
$(\rho,\mathcal{H})$ be a unitary representation of the group
$EL_k(R)$ and let $v\in \mathcal{H}$ be a unit vector s.t
$\|\rho(s)v-v\|<\epsilon$ for all $s \in S_k(R)$. Then
$\|\rho(g)v-v\|\leq 2M(l)\epsilon$ for every elementary matrix $g$.
\elem

Recently Shalom in \cite{Sh3,Sh4} gave a sufficient criterion for
$EL_d(R)$ to have Kazhdan property $(T)$.

\bthm \label {Shalom} Let $R$ be a f.g. associative ring with $1$
and with stable range $r$. Then for all $d > max{\{2,r\}}$, the
group $EL_d(R)$ has Kazhdan property $(T)$.

\ethm
 It is clear that
$M_{4n}(R)$ and  $M_4(M_n(R))$ are isomorphic as rings. Now let us
look at the multiplicative group $EL_4(M_n(R))$ contained in
$M_4(M_n(R))$. It is easy to see that each elementary matrix in
$EL_4(M_n(R))$ is a matrix in $EL_{4n}(R)$ and therefore
$EL_4(M_n(R)) \subseteq EL_{4n}(R)$. Also it is clear that
$EL_{2n}(R) \subseteq GL_2(M_n(R))$ and for $2n\geq 3$ we have
$[EL_{2n}(R),EL_{2n}(R)]=EL_{2n}(R)$. Now, since
$$[GL_2(M_n(R)),GL_2(M_n(R))] \subseteq EL_4(M_n(R))$$ (see Lemma \ref{comm} in this paper), we get
that $EL_{2n}(R) \subseteq EL_4(M_n(R))$, sitting in the upper left
block. The same is true for the case $EL_{2n}(R)$ sitting in the
lower right block. This implies that the generating set $S_{2n}(R)$
for $EL_{2n}(R)$ is a subset of $EL_4(M_n(R))$, hence
$EL_{4n}(R)=EL_4(M_n(R))$.

Let $R$ be a f.g. associative ring generated by $a_1=1,...,a_l$. Let
$$A_i=\left(\begin{array}{cccccc}
  a_i &  0  & . & . & . & 0  \\
  0 &    0  & . & . & . & 0  \\
  . & . & . & . & . & . \\
  0 & . & . & . & . & 0  \\
\end{array} \right)$$

and $$B=\left(\begin{array}{cccccc}
  0 &  1  & . & . & . & 0  \\
  0 &   0  & 1 & . & . & 0  \\
  . &   .  & . & . & . & 0  \\
  . & . & . & . & . & 1 \\
  (-1)^{n-1} & . & . & . & . & 0  \\
\end{array} \right).$$

It is easy to see that the ring $M_n(R)$ is generated by the set
$\{A_1,..,A_l,B\}$, of size $l+1$ which is independent of
$n$.

 For $n \geq 2$ and for arbitrary field $K$ it is obvious
that $SL_n(K)$ is generated by the set of elementary matrices and
thus $SL_n(K)=EL_n(K)$. Let $\mathcal{O}$ be a ring of integers in
an algebraic number field. Bass, Milnor and Serre \cite{BMS} have
shown that every matrix in $SL_n(\mathcal{O})$ for $n \geq 3$, may
be written as a product of elementary matrices, therefore
$EL_n(\mathcal{O})=SL_n(\mathcal{O})$. Moreover, Suslin in \cite
{Sus} proved that
$EL_n(\mathbb{Z}[x_1,...,x_m])=SL_n(\mathbb{Z}[x_1,...,x_m])$,
again for $n\geq 3$.

\subsection {Stable range}
\bde A sequence $\{a_1,...,a_n\}$ in a ring $R$ is said to be left
unimodular if $Ra_1+...+Ra_n=R$. In case $n \geq 2$, such a sequence
is said to be reducible if there exist $r_1,...,r_{n-1}\in R$ such
that $R(a_1+r_1a_n)+\cdots+R(a_{n-1}+r_{n-1}a_n)=R$.
 \ede

 This reduction notion leads directly to the definition of stable
 range.

 \bde A ring $R$ is said to have left stable range $\leq n$, if
 every left unimodular sequence of length $>n$ is reducible.
 The smallest such $n$ is said to be the left stable range of $R$.
 \ede

Vaserstein has proved that for any ring $R$, the left stable range
is equal to the right stable range \cite{Vas2}. Thus, we write
$sr(R)$ for this common value and call it simply the stable range
of $R$.

The reader should be aware that there is an inconsistency of $\pm
1$ in the definition of the stable range in the literature.

In \cite{Vas1} Vaserstein proves the following fundamental theorem:
 \bthm \label {Vas} Let $R$ be a ring with stable range $r$,
then the canonical mapping $$S_n:GL_n(R)/EL_n(R)\rightarrow
K_1(R)$$ is bijective for all $n \geq r +1$.
 \ethm

As a consequence of this theorem it is easy to verify that for all
$n\geq r+2$, $$EL_{n-1}(R)=EL_n(R)\bigcap GL_{n-1}(R).$$

We will need the following known fact which can be found in
[HO,4.1.18].

 \bthm \label {matrix_sr} Let $R$ be a f.g. ring with stable range $k$,
then $sr(M_n(R))=1+[\frac{k-1}{n}]$ where $[x]$ is the greatest
integer function.

 \ethm

\section{\bf Proof of Theorem \ref{main_them0} :}
We will show that every element in $EL_{4n}(R)$ is a bounded product
of elementary matrices in $EL_4(M_n(R))$ and an element of
$EL_d(R)$, where $4n >d$ and $d$ is a fixed number s.t $d$ is greater than the
stable range of the ring $R$.

In [DV, Lemma  \textbf{9}] Dennis and Vaserstein proved the
following lemma:

 \blem \label {lem_dv}
Let $R$ be an associative ring with $1$ with $sr(R)\leq r$. Then
for any $n\geq r$ we have $$GL_n(R)=ULUL\left(\begin{array}{cc}
  GL_r(R) & 0   \\
  0 & I_{n-r}   \\
\end{array} \right)=UL\left(\begin{array}{cc}
  GL_r(R) & 0   \\
  0 & I_{n-r}   \\
\end{array} \right)UL$$
where $U$ (resp., $L$) is the group of all upper (resp., lower)
triangular matrices in $EL_n(R)$ with $1$ along the main diagonal.

\elem

It is easy to see that for every $d$  with $n \geq d \geq r$ we
also have:
$$GL_n(R)=ULUL\left(\begin{array}{cc}
  GL_d(R) & 0   \\
  0 & I_{n-d}   \\
\end{array} \right)=UL\left(\begin{array}{cc}
  GL_d(R) & 0   \\
  0 & I_{n-d}   \\
\end{array} \right)UL.$$

 Here is quantitative
version which counts the number of elementary matrices for the
case $n=4$ and $sr(R)\leq 2$.

 \blem \label {elm_1} Let $R$ be an associative ring with $1$ with $sr(R)\leq 2$.
 Then every matrix in $GL_{4}(R)$ can be presented as a
 multiplication of $20$ elementary matrices and a matrix of the form
 $$\left(\begin{array}{cccc}
  A & B  & 0 & 0  \\
  C & D & 0 & 0  \\
  0 & 0 & 1 & 0 \\
  0 & 0 & 0 & 1  \\
\end{array} \right)$$ where  $$\left(\begin{array}{cc}
  A & B \\
  C & D \\
\end{array} \right) \in GL_2(R).$$ \elem

\begin{proof}

Let $M$ be an arbitrary matrix in $GL_4(R)$, so $M$ is invertible,
say
$$M=\left(\begin{array}{cccc}
  * & * & * &  a_1 \\
  ** & * & * &  a_2 \\
  ** & * & * &  a_3 \\
  ** & * & * &  a_4 \\
\end{array}
\right)$$

As $M$ is invertible, we can write
$$M^{-1}=\left(\begin{array}{cccccc}
  * & *  & * & * \\
  ** & * & * & * \\
  ** & * & * & * \\
  b_1 & b_2 & b_3 & b_4 \\
\end{array}
\right)$$

with $b_1a_1+b_2a_2+b_3a_3+b_4a_4=1$. In particular
$$Ra_1+Ra_2+Ra_3+Ra_4=R.$$
In this proof we consider only the most difficult case, where $a_i
\neq 0$ for $i=1,2,3,4$. The other cases follow immediately. Since
$sr(R)\leq 2 $, we get that there are $t_1,t_2 \in R$ such that
$$R(a_1+t_1(b_3a_3+b_4a_4))+R(a_2+t_2(b_3a_3+b_4a_4))=R.$$
Thus, there exist $x_1,x_2 \in R$ such that
$$x_1(a_1+t_1(b_3a_3+b_4a_4))+x_2(a_2+t_2(b_3a_3+b_4a_4))=1$$

For the elementary operation of adding a multiplication of a row $j$
by scalar $c$ to row $i$, we will use the following notation ${R_i
\leftarrow R_i+cR_j}$.

$$M=\left(\begin{array}{cccc}
  * & *  & * &  a_1 \\
  ** & * &  * &  a_2 \\
  ** & * &  * & a_3 \\
  ** & * &  * & a_4 \\
\end{array}
\right)  \begin{array}{cc}
  {R_1 \leftarrow R_1+t_1b_4R_4} \\
  {R_1 \leftarrow R_1+t_1b_3R_3} \\
\end{array}\left(\begin{array}{cccc}
  * & *  & *  & a_1+ t_1(b_3a_3+b_4a_4) \\
  ** & * &  * &  a_2 \\
  ** & * &  * & a_3 \\
  ** & * &  * & a_4 \\
\end{array}\right)$$
$$ \begin{array}{cc}
  {R_2 \leftarrow R_2+t_2b_4R_4} \\
  {R_2 \leftarrow R_2+t_2b_3R_3} \\
\end{array}
\left(\begin{array}{cccc}
  * & *  & * &  a_1+ t_1(b_3a_3+b_4a_4) \\
  ** & * &  * &  a_2+ t_2(b_3a_3+b_4a_4) \\
  ** & * &  * & a_3 \\
  ** & * &  * & a_4 \\
\end{array}
\right)
$$

$$  \begin{array}{cc}
  {R_4 \leftarrow R_4+(a_4-1)x_1R_1} \\
\end{array}
\left(\begin{array}{cccccc}
  * &  * & * & a_1+ t_1(b_3a_3+b_4a_4) \\
  ** & * & * &  a_2+ t_2(b_3a_3+b_4a_4) \\
  ** & * & * & a_3 \\
  ** &  * & * & a_4+(a_4-1)x_1[a_1+ t_1(b_3a_3+b_4a_4)] \\
\end{array}
\right)
$$

$$  \begin{array}{cc}
  {R_4 \leftarrow R_4+(a_4-1)x_2R_2} \\
\end{array}
\left(\begin{array}{cccc}
  * & *  & * &  a_1+ t_1(b_3a_3+b_4a_4) \\
  ** & * &  * &  a_2+ t_2(b_3a_3+b_4a_4) \\
  ** & * &  * & a_3 \\
  ** & * &  * & 1 \\
\end{array}
\right)
 $$

 This was done by $6$ elementary operations.

Now since we have $1$ in the lower right corner, we use $6$ elementary matrices to annihilate all the rest of the last row and column.  In a similar way we can use again $8$ elementary
matrices and bring $M$ to the form

$$=\left(\begin{array}{cccc}
  * & *   & 0 & 0 \\
  ** & *  & 0 &  0 \\
  0 & 0  & 1 & 0 \\
  0 & 0  & 0 & 1 \\
\end{array}
\right).$$
\end{proof}

\blem \label {lem_ULUL} Let $R$ be a f.g. associative ring with $1$
and with  $sr(R)\leq r$. If $m \geq d \geq r \geq 3$, then every
element in $GL_m(R)$ is a product of at most $8$ commutators, and an
element of the form
$$\left(\begin{array}{cc}
  S & 0   \\
  0 & I_{m-d}  \\
\end{array}
\right)$$ where $S \in GL_d(R)$. \elem

\begin{proof}
 From the remark after Lemma \ref {lem_dv} we get that any element $T\in
GL_m(R)$ can be presented as
$$T=L_1U_1L_2U_2\left(\begin{array}{cc}
  S & 0   \\
  0 & I_{m-d}  \\
\end{array}\right)$$
where $S \in GL_d(R)$. Van der Kallen \cite{KW} asserted that
every triangular matrix with entries in a ring $R$ and $1$ on the
diagonal, can be expressed as a product of three commutators. This
result was improved to two commutators by Dennis and Vaserstein [
DV, Lemma \textbf{13}].
\end{proof}

\blem \label {comm}Let $R$ be a f.g. associative ring with $1$.
Then every matrix in $GL_4(R)$ of the form
$$\left(\begin{array}{cccc}
  A & B  & 0 & 0 \\
  C & D & 0 & 0  \\
  0 & 0 & 1 & 0 \\
  0 & 0 & 0 & 1 \\
\end{array} \right)$$

where
$$\left(\begin{array}{cc}
  A & B \\
  C & D \\
\end{array} \right) \in GL_2(R),$$
is  a commutator in $GL_{2}(R)$, can be presented as a
multiplication of 40 elementary matrices in $EL_4(R)$.
 \elem

\begin{proof}
We use the following identities:

$$\left(\begin{array}{cc}
  [h_1,h_2] & 0   \\
  0 & I_{2 \times 2}  \\
\end{array}
\right) =\left(\begin{array}{cc}
  0 & h_1  \\
  -h_1^{-1} & 0 \\
\end{array}
\right) \left(\begin{array}{cc}
  0 & -h_2^{-1}    \\
  h_2 & 0  \\
\end{array}
\right) \left(\begin{array}{cc}
  (h_2h_1)^{-1} & 0   \\
  0 & h_2h_1  \\
\end{array}
\right).$$

$$\left(\begin{array}{cc}
  h & 0   \\
  0 & h^{-1}  \\
\end{array}
\right) = \left(\begin{array}{cc}
  1 & h   \\
  0 & 1  \\
\end{array}
\right) \left(\begin{array}{cc}
  1 & 0   \\
  -h^{-1} & 1  \\
\end{array}
\right)\left(\begin{array}{cc}
  1 & h-1   \\
  0 & 1  \\
\end{array}
\right)
\left(\begin{array}{cc}
  1 & 0   \\
  1 & 1  \\
\end{array}
\right)
\left(\begin{array}{cc}
  1 & -1   \\
  0 & 1  \\
\end{array}
\right)
$$

$$\left(\begin{array}{cc}
  0 & h   \\
  -h^{-1} & 0  \\
\end{array}
\right) =\left(\begin{array}{cc}
  I_{2 \times 2} & h  \\
  0 & I_{2 \times 2}  \\
\end{array}
\right) \left(\begin{array}{cc}
  I_{2 \times 2} & 0   \\
  -h^{-1} & I_{2 \times 2}  \\
\end{array}
\right) \left(\begin{array}{cc}
  I_{2 \times 2} & h   \\
  0 & I_{2 \times 2}  \\
\end{array}
\right)
$$

Now since $$\left(\begin{array}{cccc}
  1 & 0 & a & b   \\
  0 & 1 & c & d   \\
  0 & 0 &1 & 0   \\
  0 & 0 &0 & 1   \\
\end{array}
\right) = \left(\begin{array}{cccc}
  1 & 0 & 0 & 0   \\
  0 & 1 & c & 0   \\
  0 & 0 &1 & 0   \\
  0 & 0 &0 & 1   \\
\end{array}
\right) \cdot \left(\begin{array}{cccc}
  1 & 0 & 0 & 0   \\
  0 & 1 & 0 & d   \\
  0 & 0 &1 & 0   \\
  0 & 0 &0 & 1   \\
\end{array}
\right) \cdot \left(\begin{array}{cccc}
  1 & 0 & a & 0   \\
  0 & 1 & 0 & 0   \\
  0 & 0 &1 & 0   \\
  0 & 0 &0 & 1   \\
\end{array} \right)\cdot \left(\begin{array}{cccc}
  1 & 0 & 0 & b   \\
  0 & 1 & 0 & 0   \\
  0 & 0 &1 & 0   \\
  0 & 0 &0 & 1   \\
\end{array} \right)$$

and in the same way for the lower case, we get that every element
of the form $$ \left(\begin{array}{cc}
  I_{2 \times 2} & h   \\
  0 & I_{2 \times 2}  \\
\end{array}
\right)   or   \left(\begin{array}{cc}
  I_{2 \times 2} & 0  \\
  h & I_{2 \times 2}  \\
\end{array}
\right)$$

can be written as a products of at most $4$ elementary matrices
(and $2$ for the case that $h=\pm1$) in $EL_4(R)$. The result now
follows.
\end{proof}

\bpr \label{pr_el}For $n \geq d > r$, every element in
$EL_{4n}(R)$ can be presented as product of at most $340$
elementary matrices in $EL_4(M_n(R))$ and an element in $EL_d(R)$.
\epr

\brem The constant $340$ related to bounded generation can be
improved. \erem
\begin{proof}
We take an arbitrary element $T$ in $EL_{4n}(R)$. By Theorem
\ref{matrix_sr}, $sr(M_n(R))\leq 2$, hence by the argument after
Theorem \ref{Shalom} and Lemma \ref{elm_1} we get that $T$ can be
presented as a multiplication of $20$ elementary matrices of
$EL_4(M_n(R))$ and a matrix of the form
 $$\left(\begin{array}{cccccc}
  A & B  & 0 & 0  \\
  C & D & 0 & 0 \\
  0 & 0 & Id & 0  \\
  0 & 0 & 0 & Id  \\
\end{array} \right)$$ where  $$\left(\begin{array}{cc}
  A & B \\
  C & D \\
\end{array} \right) \in GL_2(M_n(R)) \subseteq GL_{2n}(R).$$

As $2n \geq sr(R) +2$, by the argument after Theorem \ref{Vas} it
follows that $$\left(\begin{array}{cc}
  A & B \\
  C & D \\
\end{array} \right) \in EL_{2n}(R).$$

 Now by Lemma \ref{lem_ULUL} (for $m=2n$) we get that
$$\left(\begin{array}{cccccc}
  A & B  & 0 & 0  \\
  C & D & 0 & 0 \\
  0 & 0 & Id & 0  \\
  0 & 0 & 0 & Id  \\
\end{array} \right).$$
 can be written as a product of $8$ commutators and an element $S \in EL_d(R)$.
 By Lemma \ref{comm}, applied to $EL_4(M_n(R)$, the result follows.
\end{proof}

We will need the following known result which is proved in many
papers (see  Ch. \textbf{3}, Cor. \textbf{11} in \cite{HV} and
Lemma \textbf{2.5} in  \cite{Sh1}):

\blem \label{cen_mas} Let $(\rho,\mathcal{H})$ be a unitary
representation of a group $\Gamma$. Suppose that for some unit
vector $v \in \mathcal{H}$, one has for all $g\in \Gamma$:
$\|\rho(g)v-v\|<\sqrt{2}$. Then there exist a non-zero
$\Gamma$-invariant vector in $\mathcal{H}$.

\elem

 Before we get to the proof of the Theorem \ref{main_them0}
we give a general property of a Kazhdan group.

\blem \label{kaz_inv} Let $\Gamma$ be a group generated by a set
$S$ with Kazhdan constant $\epsilon>0$, and let $\epsilon
> \delta >0$ be given. Let $(\rho,\mathcal{H})$ be a unitary representation of $\Gamma$ and
assume that there exists $v \in \mathcal{H}$ with $\|v\|=1$ s.t for
all $s \in S$, $\|\rho(s)v-v\| \leq \delta$. Then $v$ is
$(\Gamma,2\cdot \frac{\delta}{\epsilon})$-invariant.
 \elem

\begin{proof}

 Decompose $\rho$ into the trivial component $\sigma_0$ and the non-trivial component $\sigma_1$,
 $\rho=\sigma_0+\sigma_1$, and accordingly decompose $v=v_0+v_1$.

 For all $s \in S$,

$$
\delta \geq  \|\rho(s)v-v\|=
\|\sigma_0(s)v_0-v_0+\sigma_1(s)v_1-v_1\|=\|\sigma_1(s)v_1-v_1\|
$$

Taking the maximum over $S$, we get

\begin{equation}\label{eqt2}
\delta \geq   max _{s \in S} \|\sigma_1(s)v_1-v_1\| >
 \epsilon \cdot \|v_1\|.
\end{equation}

This implies that $\|v_1\|<\frac{\delta}{\epsilon}$,
 hence for all $g \in \Gamma$, we get $\|\rho(g)v-v\|=\|\rho(g)v_1-v_1\| < 2\cdot \frac{\delta}{\epsilon}$ as required.
\end{proof}

 Now we are ready to prove Theorem \ref{main_them0}.
Let $R$ be an associative ring generated by $l$ elements s.t.
$sr(R)\leq r$. Let $d > r$ and assume that $EL_d(R)$ is generated by
the set $F$, where $k_0=|F|$, and with Kazhdan constant $\epsilon_0$
w.r.t $F$. Notice that by Theorem \ref{Shalom} we know that
$EL_d(R)$ has Kazhdan property $(T)$.

For $EL_{4n}(R)$ where $n \geq d$ we will show that the group
$EL_{4n}(R)$ has Kazhdan constant $(k_1,\ep_1)$ with respect to the
generating set $\tilde{S}=F\bigcup S_4(M_n(R))$, where
$k_1=k_0+|S_4(M_n(R))|$ and $\ep_1=\ep_1(\ep_0,l)$ ($S_4(M_n(R))$
was defined in Subsection \ref{EL_T} and its order is independent of
$n$).

Let $(\rho, \mathcal{H})$ be a unitary representation of
$EL_{4n}(R)$ and suppose that $v \in \mathcal{H}$ is a unit vector
which is $(\tilde{S},\ep_1)$-invariant (the constant $\ep_1$ will be
determined later).

Since $EL_4(M_n(R))=EL_{4n}(R)$, we get that $v$ is
$(S_4(M_n(R)),\ep_1)$-invariant, hence by Lemma \ref{Kas_el} any
elementary matrix $g \in EL_4(M_n(R))$ satisfes:
$$\|\rho(g)v-v\|<2M(l+1)\ep_1$$
provided $\ep_1<\ep_0$.

 Now, restricting $\rho$ to the subgroup
$EL_d(R)$, we see that $v$ is $(F,\ep_1)$-invariant and hence by
Lemma \ref{kaz_inv} we get for any $g \in EL_d(R)$:
$$\|\rho(g)v-v\| < 2\cdot \frac{\ep_1}{\epsilon_0}.$$
 Let $g$ be an arbitrary element in
$EL_{4n}(R)$. By Proposition \ref{pr_el}, $g$ is expressible as a
product of at most $340$ elementary matrices of $EL_4(M_n(R))$ and
an element $g_c \in EL_d(R)$:
$$g=g_1 \cdot ... \cdot g_{340} \cdot g_c.$$

Therefore

$$\|\rho(g)v-v\| \leq
\sum_{i=1}^{340}\|\rho(g_i)v-v\|+\|\rho(g_c)v-v\| <  340\cdot
2M(l+1)\ep_1 + 2\cdot \frac{\ep_1}{\epsilon_0}.$$

If we choose  $\ep_1 < \frac{\ep_0 \cdot\sqrt{2}}{680M(l+1)+2}$ we
obtain for all $g \in EL_{4n}(R)$

$$\|\rho(g)v-v\| < \sqrt{2}.$$

Therefore $v$ is $(EL_{4n}(R),\sqrt{2})$ invariant and hence by
Lemma \ref{cen_mas} there exists $0\neq v_0 \in \mathcal{H}$ which
is $EL_{4n}(R)$-invariant.

To complete the proof, the only cases left are $EL_m(R)$ where $4n <
m < 4(n+1)$ or $ d < m<4d$. Let $m=4n+1$ (resp. $m=d+1$) and $T \in
EL_m(R)$. Since $m > sr(R)$, by using elementary operations (in the
same way we did in the proof of Lemma \ref{elm_1}) we can reduce $T$
to a matrix which has $1$ in the lower right corner. The elementary
matrices we use can be decompose to two sets:
$$\left(\begin{array}{cccccc}
  1 & 0  & \ldots & 0  & 0\\
  0 & 1 & \ldots  & 0 &  * \\
  \vdots & \vdots & \vdots & 0 & * \\
  0 & 0 & \vdots & 1 & *  \\
  0 & 0 & 0 & 0 & 1  \\
\end{array} \right) and
\left(\begin{array}{cccccc}
  1 & 0  & \ldots & 0  & *\\
  0 & 1 & \ldots  & 0 &  0 \\
  \vdots & \vdots & \vdots & 0 & 0 \\
  0 & 0 & \vdots & 1 & 0  \\
  0 & 0 & 0 & 0 & 1  \\
\end{array} \right).$$
The left set can be seen as part of the group $EL_{m-1}(R)$ which
sits in the lower right part of $EL_m(R)$. For the right set,
since all the elementary subgroups are conjugate, this set is part
of some conjugate of $EL_{m-1}(R)$ which sits (say) in the lower
right part of $EL_m(R)$. This imply that we can bring $T$ to a
matrix which has $1$ in the lower right corner by using a bounded
products of elements from groups which are isomorphic to
$EL_{m-1}(R)$, independent of $m$. Use the $1$ in the lower right
corner to annihilate all the rest of the last column and row (this
can be done is a similar way as a bounded products of groups
isomorphic to $EL_{m-1}(R)$). Now from the argument after Theorem
\ref{Vas}, it is easy to verify that $EL_m(R)$ is a bounded
product of groups isomorphic to $EL_{m-1}(R)$, independent of $m$.
The same procedure we do for  $4n< m<4(n+1)$ (resp. $d<m<4d$), and
we get that $EL_m(R)$ is a bounded product of groups isomorphic to
$EL_n(R)$ where $n=4\cdot \lfloor\frac{m}{4}\rfloor$ (reps. groups
isomorphic to $EL_d(R)$), and by Lemma \ref{BP_gp} we get uniform
Kazhdan constant for $EL_n(R)$ where $n \geq d$.

\qed

\section {Acknowledgments.} This paper is part of the author's PhD
thesis. The author wishes to thanks Moshe Jarden and Ehud de
Shalit for their comments regarding the number of generators of
ring of integers. Thanks to Noa Edelstein and Yehuda Shalom for
useful discussions. Thanks are also due to E. Bagno for his
careful reading and several helpful comments and to the anonymous referee
for his report. The author is grateful to his advisor Alex
Lubotzky for introducing him the problem and for fruitful
conversations.


Einstein Institute of Mathematics, The Hebrew University of
Jerusalem, Jerusalem 91904, Israel

E-mail address:ouzy@math.huji.ac.il

\end{document}